\documentclass[a4paper]{article}
\usepackage[utf8]{inputenc}

\usepackage[english]{babel}
\usepackage{graphicx}
\usepackage{xcolor}
\usepackage{amsmath,amsthm}
\usepackage{tikz}
\usetikzlibrary{tikzmark}
\usepackage{amssymb}
\usepackage{mathrsfs}
\usepackage{txfonts,pxfonts,tikz}
\usepackage{tgschola}
\usepackage[T1]{fontenc}
\usepackage[margin= 1 in]{geometry}
\usepackage{mathtools}
\usepackage{enumitem}   
\usepackage{listings}
\usepackage{float}

\newcommand{\defeq}{\vcentcolon=}

\usepackage[backend=biber,style=nature]{biblatex}
\usepackage{hyperref}
\hypersetup{
    colorlinks=true,
    citecolor=cyan,
    }

\newcommand{\Z}{\mathbb{Z}}
\newcommand{\N}{\mathbb{N}}
\newcommand{\F}{\mathbb{F}}

\title{An $\ell^2$ Obstruction for Elementary Embeddings of Hyperbolic Groups}
\author{Connor MacMahon}
\date{}

\theoremstyle{definition}

\theoremstyle{plain}

\begin{document}

\maketitle

\begin{abstract}
  The first $\ell^2$ Betti number of a group is non-decreasing under various embeddings arising from first order logic. Strict inequality is proved for elementary embeddings of non-abelian proper subgroups within torsion free hyperbolic groups using Perin's classification of such inclusions. The monotonicity is further demonstrated for existential embeddings of arbitrary finitely generated groups.  
\end{abstract}

\section{Introduction}

In groundbreaking work of Sela, it was demonstrated that each of the natural inclusions along the tower:
\begin{align*}
    \F_2 \hookrightarrow \F_3 \hookrightarrow \F_4 \hookrightarrow \dots
\end{align*}
consisting of free groups is elementary (\cite{FO}). This strongly resolved the Tarski problem, demonstrating that all free groups $\F_n$ for $n \geq 2$ have the same first order theory. Since this breakthrough, the structure of elementarily embedded subgroups of torsion free hyperbolic groups has been completely elucidated \cite{FO} in terms of extended towers. The purpose of this note is to use this characterization to prove that $\ell^2$ Betti numbers may arise as an obstruction to elementary embeddings.

\medskip

First defined in the seminal work of Atiyah \cite{Atiyah} in the context of index Theorems, the first $\ell^2$ Betti number of a countable discrete group has become an invariant of great interest at the intersection of algebra, geometry, and functional analysis (see \cite{CheegGrom} and \cite{Luck}). In view of the myriad interpretations that the first $\ell^2$ Betti number $\beta_1^{(2)}(G)$ of a group $G$ has, it is desirable to investigate how $\beta_1^{(2)}(G)$ is related to the first order theory of $G$. 

\medskip

The main result is as follows.

\medskip

\noindent \textbf{Theorem 1:} Suppose that $G$ is a torsion free hyperbolic group and that $H$ is a proper non-abelian subgroup which is elementarily embedded in $G$. Then $\beta_1^{(2)}(H) + \frac{1}{2} \leq \beta_1^{(2)}(G)$.

\medskip

It should be recorded that a similar homological inequality involving mod $2$ Betti numbers is known (\cite{FO} 4.5.14 and 2.0.2). 

\medskip

Upon weakening the assumptions on the inclusion, a non-strict version of the above inequality still holds. Indeed, in the case of an existential embedding one has the following.

\medskip

\noindent \textbf{Theorem 2:} Let $H$ and $G$ be finitely generated groups such that $H$ is existentially embedded within $G$. Then $\beta_1^{(2)}(H) \leq \beta_1^{(2)}(G)$. 

\medskip

The proof of Theorem 2 is substantially different from that for Theorem 1. Indeed, it uses a semi-continuity property of $\ell^2$ Betti numbers in the space of marked groups developed in \cite{sem}. 

\medskip

As an aside for readers familiar with the analogy between the first $\ell^2$ Betti number of a group and the $1$-bounded entropy of a von Neumann algebra, Theorem 2 may be viewed as exactly analogous to that in \cite{Sri} for the $1$-bounded entropy of von Neumann algebras (a notion defined in work of Jung, Shlyakhtenko, and Hayes \cite{Jung}, \cite{Reg}). Briefly, this result states that if $N \leq M$ is an existential inclusion of von Neumann algebras, i.e.
\begin{align*}
    N \leq M \leq N^\omega
\end{align*}
for a non-principal ultrafilter $\omega$ (here $N$ is embedded diagonally within $N^\omega$), then one obtains the inequality $h(N) \leq h(M)$ in the $1$-bounded entropy.

\medskip

Additional similarities between the $1$-bounded entropy $h$ and the first $\ell^2$ Betti number have been noted in various contexts. Namely, there are similar Theorems regarding their vanishing and behavior with respect to normalizers (compare \cite{PT}, \cite{Reg}, and \cite{jung2016ranktheoreml2invariantsfree}). There are even some results relating the two invariants. For instance, if a finitely presented sofic group $G$ satisfies $\beta_1^{(2)}(G)=0$, then $h(L(G))<\infty$ (\cite{Sh},\cite{Betti}). An example application of Theorem 2 regarding wq-normal subgroups is included below (see Corollary 3). 

\medskip

Now to the proofs. For preliminaries on $\ell^2$ Betti numbers, see the paper \cite{PT} and the comprehensive survey \cite{Luck}. 

\medskip

\noindent \textbf{Ackowledgements:} The author acknowledges support from NSF CAREER award DMS-21447. Moreover, the author warmly recognizes Ben Hayes and David Jekel for reading over various drafts of the paper and providing helpful suggestions.

\section{Logic Preliminaries}

Here, the necessary first order logic is introduced following \cite{ModTHy}. Within the language of groups, the terms are words in the variables denoting group elements, their inverses, and the identity element. An atomic formula consists of two terms separated by "$=$".

\medskip

To describe first order formulae, one combines atomic formulae, logical connectors ($\vee$, $\wedge$, $\neg$), and quantifiers ($\forall, \exists$). Where no free variables are present (i.e. each variable is quantified over), a formula is referred to as a sentence. For example, the following is a sentence in the first order language of groups:
\begin{align*}
    \forall g_1 \exists g_2 (g_1g_2g_1^{-1}g_2^{-1}=1) \wedge \neg(g_2 =1).
\end{align*}
The \textit{elementary theory} of a group $G$ is the collection of first order sentences which hold true in $G$. An embedding $H \hookrightarrow G$ of groups is called elementary if each first order sentence $\phi$ in the language of groups with constants in $H$ satisfies $G \models \phi$ if and only if $H \models \phi$ \cite{FO}.

\medskip

Observe that if $H$ is elementarily embedded within $G$, then $H$ and $G$ have the same elementary theory. Indeed, a first order sentence $\phi$ in the language of groups is also a first order sentence in the language of groups with constants in $H$. If $G$ and $H$ have the same elementary theory, then they are called elementarily equivalent. More precisely, $G$ and $H$ are elementarily equivalent provided that for each first order sentence $\phi$ in the language of groups, $G \models \phi$ if and only if $H \models \phi$.

\medskip

It should be noted that $H$ being elementarily embedded within $G$ is strictly stronger than $H$ being elementarily equivalent to $G$ (for instance, $2\Z \subset \Z$ is not an elementary embedding \cite{FO} 3.1.1). 

\medskip

Below, it will be desirable to also consider existential embeddings. An existential sentence is a first order sentence $\phi$ in the language of groups which takes the following form:
\begin{align*}
    \exists x_1 \exists x_2 \dots \exists x_n \psi(x_1,\dots,x_n),
\end{align*}
where $\psi$ is a quantifier free formula. The \textit{existential theory} of a group $G$ consists of the collection of existential sentences which hold true in $G$. As with elementary embeddings, an embedding $H \hookrightarrow G$ of groups is called existential provided each existential sentence $\phi$ in the language of groups with constants in $H$ satisfies $G \models \phi$ if and only if $H \models \phi$. 
 
\section{The Result}

 Before beginning, it is convenient to recall that in the present context, a \textit{graph of groups} consists of a nonempty connected graph $\Gamma$, groups $G_v$ corresponding to each vertex in $\Gamma$, and groups $G_e$ corresponding to each edge of $\Gamma$. In general, self loops at vertices are allowed, though they will not arise in the construction below. Moreover, for each edge group $G_e$, the graph of groups carries the additional data of two injective homomorphisms from $G_e$ to each of the vertex groups $G_v$ corresponding to vertices $v$ incident to the edge $e$. 

\medskip

To build the fundamental group $\pi_1(\Gamma)$ of the graph of groups $\Gamma$, fix first a maximal tree $T$ in $\Gamma$. Enumerate the edges in $T$ by $\{e_1,e_2, \dots, e_q\}$ such that each subgraph of $T$ on the edges $e_1,e_2,\dots,e_i$ for $i<q$ forms a tree. Begin with a vertex $v$ incident to $e_1$. If $v'$ is the other vertex incident to $e_1$, form the amalgamated free product $G_{v} *_{G_{e_1}} G_{v'}$ using the data of the injective edge homomorphisms $G_{e_1} \rightarrow G_v$ and $G_{e_1} \rightarrow G_{v'}$. In this way, inductively form amalgamated free products over all the edges $e_1,e_2,\dots,e_q$ in the maximal tree. Call the resulting group $G_0$.

\medskip

Contracting $\Gamma$ along $T$, one is left with a graph of groups with one vertex $v_0$ with vertex group $G_0$. The remaining edges (those not in $T$) now form loops connecting $v_0$ to itself. Inductively perform HNN extensions corresponding to each of these loops, once again using the data of the injective edge homomorphisms. The resulting group $G=\pi_1(\Gamma)$ is the fundamental group of the graph $\Gamma$ of groups.

\medskip

For more details, including a presentation of the graph of groups and its independence from the choice of maximal tree, see Serre (\cite{Trees}).

\medskip

Now, a definition which arises in Perin's characterization of elementarily embedded subgroups within torsion free hyperbolic groups.

\medskip

\noindent \textbf{Definition (\cite{FO}):} Fix a group $G$. A non-exceptional centered splitting of $G$ is a graph of groups decomposition $G = \pi_1(\Gamma)$ such that the vertices of $\Gamma$ are $v,v_1,\dots,v_r$ with $r \geq 1$, where $v$ carries the fundamental group $\pi_1(\Sigma)$ of a compact surface $\Sigma$ and each edge joins $v$ to some $v_i$.

The edge groups are isomorphic to $\Z$, and they correspond to the maximal boundary subgroups of $\pi_1(\Sigma)$. The boundary components of $\Sigma$ are in bijective correspondence with edges incident to $v$. The vertices $v_i$ are the bottom vertices, and one denotes the groups they carry by $B_i$. The vertex $v$ is called the central vertex. The abstract free product $B_1 *B_2* \dots *B_r$ is referred to as the base of the splitting. Additionally, $\chi(\Sigma) < 0$ and $\Sigma$ is not among the following list of exceptional surfaces:
\begin{itemize}
\item $S_{0,3}$ the thrice-punctured sphere,
\item $N_{1,2}$ the twice punctured projective plane,
\item $N_{2,1}$ the once punctured Klein bottle,
\item and $N_3$ the closed non-orientable surface of genus $3$. 
\end{itemize}
The closed non-orientable surface of genus $3$ is already disallowed by the fact that $r \geq 1$. However, it is still included in the list of exceptional surfaces in order to stay true to Perin's definition. 

\medskip

Topologically, one is gluing the boundary components of $\Sigma$ along loops representing nontrivial homotopy classes in the disjoint union $\bigcup_{i}X_i$ of topological spaces $X_i$ satisfying $B_i = \pi_1(X_i)$ so as to obtain a connected space. The following is a depiction of a centered splitting along with the associated graph of groups. 

\medskip

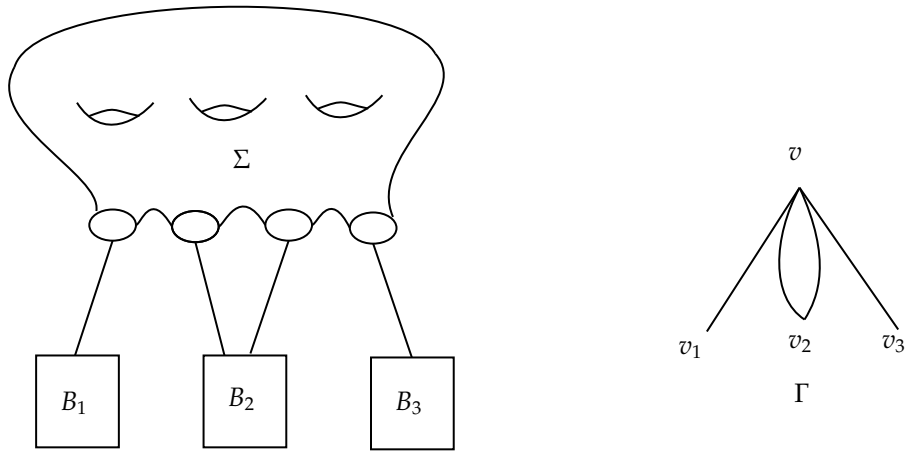
\begin{figure}[H]
\begin{center}

\tikzset{every picture/.style={line width=0.75pt}} 

\begin{tikzpicture}[x=0.75pt,y=0.75pt,yscale=-1,xscale=1]

\draw    (169.21,133.65) .. controls (172.3,115.07) and (112.25,89.8) .. (128.15,61.56) ;
\draw    (128.15,61.56) .. controls (139.63,23.66) and (320.68,21.43) .. (339.23,54.87) ;
\draw    (318.03,136.62) .. controls (304.78,106.15) and (358.66,76.43) .. (339.23,54.87) ;
\draw   (254,140.71) .. controls (254,136.4) and (259.24,132.91) .. (265.7,132.91) .. controls (272.16,132.91) and (277.4,136.4) .. (277.4,140.71) .. controls (277.4,145.02) and (272.16,148.51) .. (265.7,148.51) .. controls (259.24,148.51) and (254,145.02) .. (254,140.71) -- cycle ;
\draw   (207.19,141.45) .. controls (207.19,137.14) and (212.43,133.65) .. (218.89,133.65) .. controls (225.35,133.65) and (230.59,137.14) .. (230.59,141.45) .. controls (230.59,145.76) and (225.35,149.26) .. (218.89,149.26) .. controls (212.43,149.26) and (207.19,145.76) .. (207.19,141.45) -- cycle ;
\draw   (207.19,141.45) .. controls (207.19,137.14) and (212.43,133.65) .. (218.89,133.65) .. controls (225.35,133.65) and (230.59,137.14) .. (230.59,141.45) .. controls (230.59,145.76) and (225.35,149.26) .. (218.89,149.26) .. controls (212.43,149.26) and (207.19,145.76) .. (207.19,141.45) -- cycle ;
\draw   (165.68,140.71) .. controls (165.68,136.4) and (170.92,132.91) .. (177.38,132.91) .. controls (183.85,132.91) and (189.08,136.4) .. (189.08,140.71) .. controls (189.08,145.02) and (183.85,148.51) .. (177.38,148.51) .. controls (170.92,148.51) and (165.68,145.02) .. (165.68,140.71) -- cycle ;
\draw   (296.39,142.2) .. controls (296.39,137.89) and (301.63,134.39) .. (308.09,134.39) .. controls (314.56,134.39) and (319.8,137.89) .. (319.8,142.2) .. controls (319.8,146.51) and (314.56,150) .. (308.09,150) .. controls (301.63,150) and (296.39,146.51) .. (296.39,142.2) -- cycle ;
\draw    (189.08,140.71) .. controls (191.73,140.34) and (196.15,122.5) .. (207.19,141.45) ;
\draw    (230.59,141.45) .. controls (233.24,141.08) and (241.19,119.53) .. (254,140.71) ;
\draw    (277.4,140.71) .. controls (280.05,140.34) and (286.23,121.76) .. (296.39,142.2) ;
\draw    (159.5,79.4) .. controls (172.3,99.46) and (190.85,87.57) .. (197.92,79.4) ;
\draw    (165.68,86.09) .. controls (179.37,80.88) and (176.72,83.11) .. (190.85,86.09) ;
\draw    (216.02,77.17) .. controls (228.83,97.23) and (247.37,85.34) .. (254.44,77.17) ;
\draw    (222.2,83.86) .. controls (235.89,78.66) and (233.24,80.88) .. (247.37,83.86) ;
\draw    (274.31,75.68) .. controls (287.12,95.75) and (305.66,83.86) .. (312.73,75.68) ;
\draw    (280.49,82.37) .. controls (294.18,77.17) and (291.53,79.4) .. (305.66,82.37) ;
\draw   (139,206) -- (180.5,206) -- (180.5,252) -- (139,252) -- cycle ;
\draw   (223,206) -- (264.5,206) -- (264.5,252) -- (223,252) -- cycle ;
\draw   (307,207) -- (348.5,207) -- (348.5,253) -- (307,253) -- cycle ;
\draw    (177.38,148.51) -- (158.5,206) ;
\draw    (218.89,149.26) -- (233.5,206) ;
\draw    (265.7,148.51) -- (246.5,205) ;
\draw    (308.09,150) -- (327.5,207) ;
\draw    (522,122) -- (475.5,194) ;
\draw    (522,122) -- (571.5,193) ;
\draw    (522,122) .. controls (508.5,147) and (507.5,177) .. (524.5,188) ;
\draw    (522,122) .. controls (538.5,158) and (531.5,177) .. (524.5,188) ;

\draw (236,101.4) node [anchor=north west][inner sep=0.75pt]    {$\Sigma $};
\draw (150,222.4) node [anchor=north west][inner sep=0.75pt]    {$B_{1}$};
\draw (234,220.4) node [anchor=north west][inner sep=0.75pt]    {$B_{2}$};
\draw (318,222.4) node [anchor=north west][inner sep=0.75pt]    {$B_{3}$};
\draw (518,218.4) node [anchor=north west][inner sep=0.75pt]    {$\Gamma $};
\draw (515,100.4) node [anchor=north west][inner sep=0.75pt]    {$v$};
\draw (460,197.4) node [anchor=north west][inner sep=0.75pt]    {$v_{1}$};
\draw (515,194.4) node [anchor=north west][inner sep=0.75pt]    {$v_{2}$};
\draw (562,194.4) node [anchor=north west][inner sep=0.75pt]    {$v_{3}$};

\end{tikzpicture}

\caption{A centered splitting with surface $\Sigma$, bottom groups $B_1,B_2,B_3$, and with the associated graph of groups $\Gamma$ depicted. Reproduced from \cite{FO} Figure 4.1.}
\end{center}
\end{figure}

Additionally, the splitting $\Gamma$ is assumed to be minimal in the sense that if $v_i$ belongs to only one edge, then the edge group $\Z$ is a proper subgroup of $B_i$. Moreover, redundant vertices are permitted: cyclic vertices $v_i$ of valence $2$ with both incident edge groups equal to the vertex group are allowed.

\medskip

The importance of this construction for the present purpose is the following result, which follows from \cite{FO} 2.0.2 and the relevant definitions.

\medskip

\noindent \textbf{Theorem 3:} Suppose that $G$ is a torsion free hyperbolic group and that $H$ is a non-abelian subgroup. If $H$ is elementarily embedded in $G$, then there exists a finite chain of groups:
\begin{align*}
    G=G_0 > G_1 > \dots > G_n=H
\end{align*}
such that for each $0 \leq k \leq n-1$, either $G_{k} \cong G_{k+1} *\Z$ or $G_k$ admits a non-exceptional centered splitting with base isomorphic to $G_{k+1}$. Additionally, whenever the latter case arises one has that the genus $g$ of the surface $\Sigma$ corresponding to the central vertex in the splitting satisfies $g \geq n_1$, where $n_1$ denotes the number of bottom vertices of valence $1$ (\cite{FO} Remark 4.5.5).

\medskip

\noindent \textbf{Remark:} The existence of such a chain is necessary for an elementary embedding of a non-abelian subgroup within a torsion free hyperbolic group, but far from sufficient. In \cite{FO} section 4.5 the additional hypotheses defining an \textit{extended tower} are mentioned, upon which the claim becomes an if and only if statement. However, these extra hypotheses will not be of use here. 

\medskip

Now, the necessary results regarding $\ell^2$ Betti numbers are collected. Firstly, the following result which is proved in the appendix to \cite{FPL} in the case $n=2$ will be useful. The general case stated here follows by induction.

\medskip

\noindent \textbf{Lemma 1 (\cite{FPL} A.1):} Suppose that $G_1,G_2,\dots,G_n$ are infinite groups. Then:
\begin{align*}
    \beta_1^{(2)}(G_1* G_2 * \dots*G_n) = (n-1) + \sum_{i=1}^n \beta_1^{(2)}(G_i).
\end{align*}

\medskip

Furthermore, a Theorem of Peterson and Thom which allows for lower estimates on the first $\ell^2$ Betti number of a group $G$ from a presentation of $G$ will be relevant below. It is stated here for the convenience of the reader.

\medskip

\noindent \textbf{Theorem 4 (\cite{PT} 3.2):} Let $G$ be an infinite countable discrete group. Assume that there exist subgroups $G_1,\dots,G_n$ such that:
\begin{align*}
    G=\langle G_1,G_2,\dots,G_n \, | \, r_1^{w_1},\dots,r_k^{w_k} \rangle
\end{align*}
for elements $r_1,\dots,r_k \in G_1*G_2*\dots*G_n$ and positive integers $w_1,\dots,w_k$. Suppose that the presentation is irredundant in the sense that $r_i^l \neq 1 \in G$ for all $1<l<w_i$ and $1 \leq i \leq k$. Then the following inequality holds:
\begin{align*}
    \beta_1^{(2)}(G) \geq n-1 + \sum_{i=1}^n \bigg(\beta_1^{(2)}(G_i) - \frac{1}{|G_i|}\bigg) - \sum_{j=1}^k \frac{1}{w_j}.
\end{align*}
By convention, if $G_i$ is infinite then $1/|G_i|$ is taken to be $0$. 
\medskip

\noindent \textbf{Theorem 1:} Suppose that $G$ is a torsion free hyperbolic group and that $H$ is a proper non-abelian subgroup which is elementarily embedded in $G$. Then $\beta_1^{(2)}(H) + \frac{1}{2} \leq \beta_1^{(2)}(G)$.
\begin{proof}
    By the Theorem 3 and induction, it suffices to treat 2 cases. Firstly, if $G \cong H *\Z$, then one has $\beta_1^{(2)}(G) = \beta_1^{(2)}(H)+1$ by Lemma 1. Indeed, $H$ is infinite as $G$ is torsion free. 

    \medskip
    
    Now, suppose instead that $G$ admits a non-exceptional centered splitting $G=\pi_1(\Gamma)$ with base $B_1*B_2* \dots * B_r \cong H$ and central group $\pi_1(\Sigma)$. Moreover, the genus $g$ of $\Sigma$ satisfies $g \geq n_1$, with $n_1$ the number of bottom vertices of valence $1$.

    \medskip
    
    Suppose that $\Sigma$ has $b$ boundary components and observe that connectedness of $\Gamma$ requires that $b \geq r$. In the associated graph of groups, a maximal tree $T$ consists of exactly one edge from $v$ attached to each $v_i$. Hence, there are $r$ edges in the maximal tree, and there are $b-r$ edges not in the maximal tree. Indeed, the graph of groups has $b$ edges: one for the gluing of each boundary component.

    \medskip
    
    Let $c_1,\dots, c_b$ denote the generators of the maximal boundary subgroups of $\pi_1(\Sigma)$, ordered such that the first $b-r$ of them have their associated edges not included in the maximal tree $T$. Now, by Bass-Serre theory $G$ may be presented as follows (\cite{Trees} 5.1(b), the presentation here is shorter than Serre's as Serre insists on a symmetric generating set):
    \begin{align*}
        G = \langle \pi_1(\Sigma),B_1,B_2,\dots,B_r,t_{1},t_{2},\dots,t_{b-r}  \, | \, t_i c_i t_i^{-1}a_i^{-1},c_ja_j^{-1}, 1 \leq i \leq b-r,b-r < j \leq b \rangle.
    \end{align*}
    Here, $a_i$ is a member of the bottom group connected to the edge corresponding to the boundary component related to $c_i$ for each $i$ (same for $j$). Now, by Theorem 4 \cite{PT} with $w_i=1$ for all $i$:
    \begin{align*}
        \beta_1^{(2)}(G) \geq [r+1+(b-r)]-1 + \beta_1^{(2)}(\pi_1(\Sigma)) +  \sum_{\ell=1}^r \beta_1^{(2)}(B_\ell) -b  = \beta_1^{(2)}(\pi_1(\Sigma)) +  \sum_{\ell=1}^r \beta_1^{(2)}(B_\ell).
    \end{align*}
    Indeed, the first $\ell^2$ Betti number of a cyclic group vanishes and all involved groups are infinite as $G$ is torsion free. Recall that $H \cong B_1*B_2*\dots *B_r$, the base of the splitting, so that applying Lemma 1 yields:
    \begin{align*}
         \beta_1^{(2)}(G) \geq \beta_1^{(2)}(\pi_1(\Sigma)) + \beta_1^{(2)}(H) - (r-1).
    \end{align*}
    Employing again Theorem 4 with a standard presentation of $\pi_1(\Sigma)$ one obtains $\beta_1^{(2)}(\pi_1(\Sigma)) \geq 2g+b-2$ or  $\beta_1^{(2)}(\pi_1(\Sigma)) \geq g+b-2$ according as $\Sigma$ is orientable or not. That is, $\beta_1^{(2)}(\pi_1(\Sigma)) \geq -\chi(\Sigma)$. Observe that the least of the two possibilities is $g+b-2$. Therefore,
    \begin{align*}
        \beta_1^{(2)}(G) \geq (g+b-2) + \beta_1^{(2)}(H) - (r-1) = (g+b-r-1) + \beta_1^{(2)}(H). 
    \end{align*}
    Now, recall that $g \geq n_1$, where $n_1$ is the number of bottom vertices of valence $1$. Observe that the number of bottom vertices not of valence $1$ is then $r-n_1$. Each such vertex has valence at least $2$, while the remaining $n_1$ bottom vertices have valence $1$. Therefore:
    \begin{align*}
        2(r-n_1)+n_1 \leq b,
    \end{align*}
    the total number of edges. Rearranging yields $b+n_1 \geq 2r$. Therefore, $b+g \geq 2r$ as $g \geq n_1$. Substituting this into the above gives:
    \begin{align*}
        \beta_1^{(2)}(G) \geq (r-1)+\beta_1^{(2)}(H).
    \end{align*}
    Hence, $r-1 \leq \beta_1^{(2)}(G) - \beta_1^{(2)}(H)$. Substituting this into the result of Lemma 1 above yields:
    \begin{align*}
        \beta_1^{(2)}(G) \geq \beta_1^{(2)}(\pi_1(\Sigma)) + \beta_1^{(2)}(H) - (r-1) \geq -\chi(\Sigma) + \beta_1^{(2)}(H) - (\beta_1^{(2)}(G) - \beta_1^{(2)}(H)).
    \end{align*}
    Rearranging, 
    \begin{align*}
        \beta_1^{(2)}(G) \geq -\frac{1}{2}\chi(\Sigma) + \beta_1^{(2)}(H).
    \end{align*}
    By the definition of a centered splitting, $\chi(\Sigma)<0$. Since $\chi(\Sigma)$ is an integer, the claim is proved. 
\end{proof}

\medskip

\noindent \textbf{Example:} In this example, the above argument is worked out for the non-exceptional centered splitting in Figure $1$. Choose the maximal tree $T$ as follows.

\begin{figure}[H]
        
\begin{center}
\tikzset{every picture/.style={line width=0.75pt}} 

\begin{tikzpicture}[x=0.75pt,y=0.75pt,yscale=-1,xscale=1]

\draw [color={rgb, 255:red, 74; green, 144; blue, 226 }  ,draw opacity=1 ]   (309,96) -- (262.5,168) ;
\draw [color={rgb, 255:red, 74; green, 144; blue, 226 }  ,draw opacity=1 ]   (309,96) -- (358.5,167) ;
\draw [color={rgb, 255:red, 74; green, 144; blue, 226 }  ,draw opacity=1 ]   (309,96) .. controls (295.5,121) and (294.5,151) .. (311.5,162) ;
\draw    (309,96) .. controls (325.5,132) and (318.5,151) .. (311.5,162) ;

\draw (270,197.4) node [anchor=north west][inner sep=0.75pt]    {$\Gamma $};
\draw (302,74.4) node [anchor=north west][inner sep=0.75pt]    {$v$};
\draw (247,171.4) node [anchor=north west][inner sep=0.75pt]    {$v_{1}$};
\draw (302,168.4) node [anchor=north west][inner sep=0.75pt]    {$v_{2}$};
\draw (349,168.4) node [anchor=north west][inner sep=0.75pt]    {$v_{3}$};
\draw (327,198.4) node [anchor=north west][inner sep=0.75pt]    {$\textcolor[rgb]{0.29,0.56,0.89}{T}$};

\end{tikzpicture}
\end{center}
\caption{A maximal tree $T$ within the graph of groups $\Gamma$ from Figure 1 is colored in blue.}
    \end{figure}
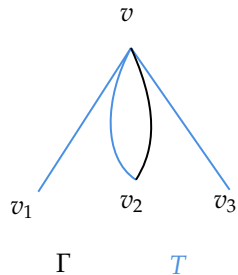
Let $c_1,c_2,c_3,c_4$ be generators for the maximal boundary subgroups of $\Sigma$, going from left to right in Figure $1$ (this is a different labelling scheme than that which is used in the above proof). Now, $G$ may be presented as follows in this case:
\begin{align*}
    G = \langle \pi_1(\Sigma),B_1,B_2,B_3,t \, | \, c_1a_1^{-1},c_2a_2^{-1},tc_3t^{-1}a_3^{-1},c_4a_4^{-1}  \rangle.
\end{align*}
Indeed, the black edge $\{v,v_2\}$ is the only edge not contained in the maximal tree $T$. Notice that there are $5$ groups generating the presentation and there are $4$ relations (the number of boundary components). Moreover, each of the involved groups is infinite. Hence, Theorem 4 \cite{PT} and Lemma 1 yield:
\begin{align*}
    \beta_1^{(2)}(G) \geq 5-1 + \beta_1^{(2)}(\pi_1(\Sigma)) + \beta_1^{(2)}(B_1)+\beta_1^{(2)}(B_2)+\beta_1^{(2)}(B_3) + \beta_1^{(2)}(\langle t \rangle) - 4 \\  = \beta_1^{(2)}(\pi_1(\Sigma)) + \beta_1^{(2)}(H) - 2,
\end{align*}
where $H$ is isomorphic to the base $B_1*B_2*B_3$ of the splitting. Now, the surface $\Sigma$ in Figure $1$ is orientable with genus $3$ and $4$ boundary components, so one obtains:
\begin{align*}
     \beta_1^{(2)}(G) \geq (2(3)+4-2) + \beta_1^{(2)}(H) -2 = 6+\beta_1^{(2)}(H). 
\end{align*}
This completes the computation.

\section{Existential Embeddings}

In this section, a proof of Theorem 2 is supplied. For preliminaries on the space of marked groups, see \cite{Marked}.

\medskip

\noindent \textbf{Theorem 2:} Let $H$ and $G$ be finitely generated groups such that $H$ is existentially embedded within $G$. Then $\beta_1^{(2)}(H) \leq \beta_1^{(2)}(G)$. 
\begin{proof}
    Denote by $V_1 = \{h_1,h_2,\dots,h_m\}$ a finite set of generators for $H$, and by $V_2 = \{g_1,g_2,\dots,g_n\}$ a finite set of generators for $G$. Set $S= V_1 \cup V_2$, a new generating set for $G$. Fix $r>0$ and consider the existential sentence:
    \begin{align*}
        \exists g_1^{(r)}\exists g_2^{(r)} \exists \dots \exists g_n^{(r)}  \bigwedge_{|w| \leq r} (w(h_1,h_2,\dots,h_m,g_1^{(r)},g_2^{(r)},\dots,g_n^{(r)})=1 \leftrightarrow w(h_1,h_2,\dots,h_m,g_1,g_2,\dots,g_n)=1).
    \end{align*}
   The conjunction runs over words $w$ in $\F(S)$ of length not exceeding $r$. This is a slight abuse of notation, as $g_1,g_2,\dots,g_n$ are not necessarily constants in the signature of the language of groups with constants in $H$. However, the formula is understood to be the conjunction of formulas of the form:
    \begin{align*}
       & w(h_1,h_2,\dots,h_m,g_1^{(r)},g_2^{(r)},\dots,g_n^{(r)})=1 \, \, \, \text{or} \\
       & \neg (w(h_1,h_2,\dots,h_m,g_1^{(r)},g_2^{(r)},\dots,g_n^{(r)})=1)
    \end{align*}
    according as $w(h_1,h_2,\dots,h_m,g_1,g_2,\dots,g_n)=1$ is true in $G$ or not. Since the embedding $H \leq G$ is existential, the sentence holds true in $H$ such that there exist members $g_1^{(r)},g_2^{(r)},\dots,g_n^{(r)}$ of $H$ satisfying $w(h_1,h_2,\dots,h_m,g_1^{(r)},g_2^{(r)},\dots,g_n^{(r)})=1$ in $H$ if and only if $w(h_1,h_2,\dots,h_m,g_1,g_2,\dots,g_n)=1$ in $G$. Denote $V_2^{(r)} \defeq \{g_1^{(r)},g_2^{(r)},\dots,g_n^{(r)}\}$, and $S^{(r)} \defeq V_1 \cup V_2^{(r)}$, a generating set for $H$. By definition,
    \begin{align*}
        \lim_{r \rightarrow \infty} (H,S^{(r)}) = (G,S)
    \end{align*}
    in the space of marked groups. In other words, $(H,S^{(r)})$ has the same relations of length at most $r$ as $(G,S)$, such that $(H,S^{(r)})$ is at distance at most $e^{-r}$ from $(G,S)$ in the space of marked groups. Hence, applying the semi-continuity of the first $\ell^2$ Betti number in the space of marked groups (\cite{sem}), the claim is proved. 
\end{proof}

An alternative proof involving ultrapowers is as follows.

\begin{proof}
    By the ultrapower characterization of existential embeddings, there exists a non-principal ultrafilter $\omega$ an an embedding $\iota: G \hookrightarrow H^\omega$ such that the image $\iota(G)$ contains the diagonal copy of $H$ within $H^\omega$. Henceforth, identify $H$ with its diagonal copy within $H^\omega$ and $G$ with its image in $H^\omega$. Let $V_1$ be a finite set of generators for $H$ and $V_2$ a finite set of generators for $G$. Denote $S=V_1 \cup V_2$, which remains a finite set of generators for $G$. 

    \medskip

    Following the proof in \cite{Marked} 6.4, write each member $s$ of $S$ as a sequence $(s_k)_{k \in \N}$. That is, lift each generator in $S$ to the product $H^\N$. For each $k$, set $S_k \defeq \{s_k \, | \, s \in S\}$. Since $S$ contains the set $V_1$ of generators for the diagonal copy of $H$, one has that for each $k$, the set $S_k$ forms a generating set for $H$. Following the proof of \cite{Marked} 6.4, there exists a subsequence $(k_i)_{i \in \N}$ such that:
    \begin{align*}
       \lim_{i \rightarrow \infty} (H,S_{k_i}) =  (G,S)
    \end{align*}
    in the space of of marked groups. Hence, applying the semi-continuity of the first $\ell^2$ Betti number in the space of marked groups (\cite{sem}), the claim is proved. 
\end{proof}

\section{Discussion and Outlook}

First, some discussion of Theorem 1. Consider its contrapositive. This states that if $H$ is a non-abelian proper subgroup of a torsion free hyperbolic group $G$, and $\beta_1^{(2)}(H)+ \frac{1}{2} > \beta_1^{(2)}(G)$, then $H$ is not elementarily embedded within $G$.

\medskip

More explicitly, there exists a first order sentence $\phi$ with constants in $H$ such that $\phi$ holds true in $G$ but not in $H$ (if instead there were a sentence which held true in $H$ but not in $G$, its negation would do). In the special case in which $\beta_1^{(2)}(G)=0$, one obtains the following Corollary.

\medskip

\noindent \textbf{Corollary 1:} Suppose $G$ is a torsion free hyperbolic group with $\beta_1^{(2)}(G) = 0$. Then $G$ admits no proper elementarily embedded subgroup $H$ which is non-abelian.

\medskip

Some torsion free hyperbolic groups $G$ with $\beta_1^{(2)}(G)=0$ are those with property (T) \cite{Betti} and the fundamental groups of closed, hyperbolic $3$-manifolds (\cite{Dod}). The first $\ell^2$ Betti number of $G$ also vanishes if $G$ has a finitely generated infinite normal subgroup of infinite index \cite{PT}. 

\medskip

Furthermore, examining the proof of Theorem 1, note that each step in the extended tower of Perin adds at least $1/2$ to the first $\ell^2$ Betti number. A straightforward induction yields the following Corollary.

\medskip

\noindent \textbf{Corollary 2:} Suppose that $H$ is a non-abelian subgroup of a torsion free hyperbolic group $G$ which is elementarily embedded within $G$. Then $G$ is an extended tower over $H$ of length at most
\begin{align*}
    2(\beta_1^{(2)}(G)-\beta_1^{(2)}(H)).
\end{align*}

\medskip

Turning now to Theorem 2, one may also consider the logical implications of its contrapositive. If $G$ is a finitely generated group and $H$ a finitely generated subgroup of $G$ with $\beta_1^{(2)}(H) > \beta_1^{(2)}(G)$, then $H$ is not existentially embedded within $G$.

\medskip

More explicitly, there exists an existential sentence $\phi$ with constants in $H$ such that $\phi$ holds true in $G$ but not in $H$. An immediate Corollary is that a free group $\F_n$ for $n \geq 2$ cannot be existentially embedded within a group $G$ with $\beta_1^{(2)}(G)=0$.

\medskip

As a final application, a Corollary which is similar in nature to the results in \cite{Reg} regarding normalizers of von Neumann algebras within ultrapowers is derived. First, a definition.

\medskip

\textbf{Definition:} A subgroup $H$ of a group $G$ is called q-normal provided there exists a generating set $S$ of $G$ such that $gHg^{-1} \cap H$ is infinite for each member $g$ of $S$. A subgroup $H$ of a group $G$ is called wq-normal if there exists an ordinal number $\alpha$, and an ascending $\alpha$-chain of subgroups such that $H_0=H$, $H_\alpha = G$, and 
\begin{align*}
    \bigcup_{\beta < \delta} H_\beta \subset H_\delta
\end{align*}
is q-normal for each ordinal $\delta$. 

\medskip

\noindent \textbf{Corollary 3:} Suppose that $H \leq G$ is an existential inclusion of finitely generated groups. If $K$ is an infinite wq-normal subgroup of $G$, then $\beta_1^{(2)}(H) \leq \beta_1^{(2)}(K)$.
\begin{proof}
    Apply \cite{PT} Theorem 5.6 and Theorem 2. 
\end{proof}

This is analogous to the claim in \cite{Reg} about normalizers of II$_1$ factors within ultrapowers due to the ultrapower characterization of existential embeddings.

\printbibliography

\end{document}